\documentclass{article}

\usepackage{amssymb,dsfont,latexsym}

\newtheorem{thm}{Theorem}
\newtheorem{defin}[thm]{Definition}
\newtheorem{prop}[thm]{Proposition}
\newtheorem{lem}[thm]{Lemma}
\newtheorem{cor}[thm]{Corollary}
\newtheorem{ex}[thm]{Example}

\newcommand\enu[1]{\smallskip\newline\makebox[5mm][l]{\rm(#1)}}

\newcommand\bp{\noindent{\it Proof.}\ }

\newcommand\2{{\frac{1}{2}}}

\begin{document}

\author{Erling St{\o}rmer}

\date{08-13-2013}

\title{A decomposition theorem for positive maps, and the projection onto a spin factor}

\maketitle 

It is shown that each positive map between matrix algebras is the sum of a maximal decomposable map and an atomic map which is both optimal and co-optimal. The result is studied  in detail for the projection onto a spin factor.

\section*{ Introduction}

The structure of positive maps between $C^*$-algebras, even in the finite dimensional case,  is still poorly understood.  The only maps which are well understood are the decomposable ones , which are sums of completely positive  and co-positive maps, hence in the finite dimensional case, are sums of maps of the form $Adv$ and $t \circ Adv$, where $t$ is the transpose map, and $Adv$ the map $x\to v^*xv$.  In the present paper we shall shed some light on the structure of positive maps by showing that they are the sum of a maximal decomposable map and an atomic map, which is bi-optimal, i.e. it majorizes neither a non-zero completely positive map nor a co-positive map. 
 
  In order to obtain a deeper understanding of this decomposition we study it in detail in Section 2 for the trace invariant positive projection of the full matrix algebra $M_{2^n}$ onto a spin factor inside it.  We shall obtain explicit formulas for the decomposable map and the bi-optimal map in the decomposition when the spin factor is irreducible and contained in the $2^{n-1} \times 2^{n-1} $ matrices over the quaternions.
  
  For the reader«s convenience we recall the main definitions concerning positive maps, see also \cite{S}.  We let $A$ be a finite dimensional $C^*$-algebra and $B(H)$ the bounded operators on a finite dimensional Hilbert space $H$.
  
  Let $\phi\colon A\to B(H)$ be a linear map. $\phi$ is \textit{positive}, written $\phi \geq 0$ or $0 \leq \phi$ if it carries positive operators to positive operators.  If $\psi$ is another positive map, $\psi$ \textit{majorizes }$\phi$ , written $\psi\geq\phi$ if $\psi - \phi \geq 0$. $\phi$ is \textit{k-positive} if $\iota_k \otimes \phi \colon {M_k \otimes A }\to {M_k \otimes B(H)}$ is positive, where $\iota_k$ is the identity map on the $k\times k$ matrices $M_k$.  $\phi$ is \textit{completely positive} if $\phi$ is k-positive for all k.  Let $t$ denote the transpose map on $B(H)$ with respect to some fixed orthonormal basis.  Then $\phi$ is \textit{k-co-positive}, (resp. \textit{co-positive})  if $t\circ\phi$ is k-positive (resp. completely positive). $\phi$ is \textit{k-decomposable} (resp. \textit{decomposable}) if $\phi$ is the sum of a k-positive and a k-co-positive map (resp. completely positive and a co-positive map).  $\phi$ is \textit{atomic} if $\phi$ is not 2-decomposable.  $\phi$ is \textit{extremal} or just extreme, if $\phi\geq\psi$ for a positive map $\psi$ implies $\psi = \lambda \phi$ for some nonnegative number $\lambda$. $\phi$ is \textit{optimal} (resp. \textit{co-optimal}) if $\phi\geq\psi$ for $\psi$ completely positive (reps. co-positive) implies $\psi =0$.  Combining the last two concepts we introduce the following definition, which has also been introduced by Ha and Kye \cite{HK}.

\begin{defin}
 $\phi$ is \textit{bi-optimal} if $\phi$ is both optimal and co-optimal.
\end{defin}

The author is grateful to E. Alfsen for many helpful discussions on spin factors.

\section { The decomposition theorem}

Let $K$ and $H$ be finite dimensional Hilbert spaces.  In \cite {M}, Theorem 3.4 Marciniak showed the surprising result that if $\phi$ is a 2-positive map (resps. 2-co-positive) which is extremal, then $\phi$ is completely positive (reps. co-positive).  His proof, see also \cite {S}, Theorem 3.3.7, contained more information, namely the following result.

\begin{lem}\label{lem1}
Let $\phi$ be a non-zero 2-positive map of $B(K)$ into $B(H)$.  Then there exists a non-zero completely positive map $\psi\colon B(K) \to B(H)$ such that $\phi\geq\psi$.
\end{lem}

A slight extension of the above lemma yields the following.

\begin{prop}\label{prop2}
Let $A$ be a finite dimensional $C^*$-algebra and $\phi\colon A \to B(H)$ a non-zero 2-decomposable map.  Then there exists a non-zero decomposable map  
$\psi\colon A \to B(H)$ such that $\phi\geq\psi$.
\end{prop}
\bp
We first consider the case when $A=B(K)$.  Since $\phi$ is 2-decomposable there exist a 2-positive map $\phi_1$ and a 2-co-positive map $\phi_2$ such that $\phi=Ê\phi_1 + \phi_2$. By Lemma 2 there is a completely positive map $\psi_1$, non-zero if $\phi_1$ is non-zero, such the $\phi_1 \geqÊ\psi_1$.  Applying Lemma 2 to $t\circ \phi_2$ we find a co-positive map $\psi_2\leq \phi_2$.  Thus $\phi\geq \psi_1 + \psi_2$,  proving the proposition when $A=B(K)$.

In the general case let $e_1,...,e_m$ be the minimal central projections in $A$, so $A= \bigoplus {_i}^m Ae_i$.  Then each $Ae_i$ is isomorphic to some $B(K)$, and $\phi_{ \vert Ae_i}$ is 2-decomposable.  By the first part $\phi_{\vert Ae_i} \geq \alpha_i + \beta_i$ with $\alpha_i$ completely positive and $\beta_i$ co-positive.  Let $\alpha = \sum \alpha_i$ and $\beta = \sum \beta_i$. Then $\alpha$ is completely positive and $\beta$ co-positive, hence $\alpha + \beta$ is a decomposable map majorized by $\phi$, completing the proof of the proposition.

\begin{cor}\label{cor3}
Each bi-optimal map of a finite dimensional $C^*$-algebra into $B(H)$ is atomic.
\end{cor}
\bp
By definition a map $\phi$ is atomic if it is not 2-decomposable.  By definition of being bi-optimal such a map $\phi$ cannot majorize a decomposable map, hence by Proposition 3 $\phi$ cannot be 2-decomposable, completing the proof.

\medskip

Since completely positive maps are sums of maps of the form $Adv$, and each co-positive map a sum of maps $t\circ Adv$, our next result reduces much of the study of positive maps to that of bi-optimal maps.  If $\phi\colon A \to B(H)$ is positive, $A$ a $C^*$-algebra, we say a decomposable map $\alpha\colon A \to B(H), \alpha\leq \phi$ is a \textit{maximal decomposable map majorized by} $\phi$ if there is no decomposable map $\psi\colon A \to B(H)$ such that $\psi \neq \alpha$, and $\alpha\leq\psi\leq\phi$.

\begin{thm}\label{thm4}
Let $A$ be a finite dimensional $C^*$-algebra and $H$ a finite dimensional Hilbert space. Let $\phi\colon A\to B(H)$ be a positive map.  Then there are a maximal decomposable map $\alpha\colon A \to B(H)$ majorized by $\phi$ and a bi-optimal, hence atomic, map $\beta\colon A \to B(H)$ such that $\phi = \alpha + \beta$.
\end{thm}
\bp
We first assume $A=B(K)$ for a finite dimensional Hilbert space $K$. Let
$$
C= \{ \psi\colon B(K) \to B(H): \psi \ \ decomposable, \psi\leq \phi\}.
$$
Then $C$ is bounded and norm closed, hence is compact in the norm topology, as $K$ and $H$ are finite dimensional.  Furthermore $C$ is an ordered set with the usual ordering on positive maps.  We show $C$ has a maximal element.  For this let $X = \{ \phi_v \in C: v \in F\}$ be a totally ordered set with $\phi_v \leq \phi_{v^{\prime}}$ if $v \leq v^{\prime}$ in $F$.  For each $v \in F$ let $X_v = \{\phi_{v^{\prime}} \in X: v \leq v^{\prime}\}$. Then $X_v$ is closed, and $X_v \supset X_{v^{\prime}}$ if $v \leq v^{\prime}$. Since $X$ is totally ordered it follows that the sets $X_v$ with $v \in F$ have the finite intersection property. Thus the intersection $\bigcap_{vÊ\in F}  X_v \neq \emptyset$, hence a map $\psi \in \bigcap X_v$ is an upper bound for  $X$.  By Zorn's lemma $C$ has a maximal element $\alpha$.  Since $C$ is closed, $\alpha$ is decomposable, $\alpha \leq \phi$, and there is no decomposeosable map $\psi\colon B(K) \to B(H)$ different from $\alpha$ such that $\alpha \leq \psi \leq \phi$.  Thus $\alpha $ is maximal decomposable map majorized by $\phi$.

Let $\beta = \phi -\alpha$.  Then $\beta$ is bi-optimal, for if $\gamma \leq \beta, \gamma \neq 0$ and decomposable, then $\alpha + \gamma$ is decomposable, and $\alpha + \gamma \leq \alpha + \beta = \phi$, contradicting maximality of $\alpha$.  Thus $\gamma = 0$, and $\beta$ is bi-optimal.

In the general case we imitate the proof of Proposition 3 and write $A$ as $A = \bigoplus Ae_i$ where the $e_i$ are minimal central projections in $A$, so $Ae_i$ is isomorphic to some $B(K)$, and we apply the first part of the proof to each $Ae_i$ in the same way as we did in the proof of Proposition 3.  The proof is complete.

\medskip

If we do not require $\alpha$ in the theorem to be maximal decomposable we can have different decompositions.  For example, if $\phi$ is a bi-optimal map, and $Tr$ is the trace on $B(K)$ , then the map $\psi(x) = \phi(1)Tr(x) + \phi(x)$ is super-positive, hence in particular completely positive, see \cite {S}, Theorem 7.5.4.  But $\psi$ has a decomposition $\psi = \alpha + \beta$, where $\alpha = \phi(1)Tr$ is completely positive, and $\beta = \phi$ is bi-optimal.

\begin{cor}\label{cor 5}
With assumptions as in Theorem 5, if $\phi$ is extreme, then $\phi$ is either of the form $Adv, t\circ Adv$ or $\phi $ is bi-optimal, so atomic.
\end{cor}

If we in the proof of Theorem 5 replace decomposable map by completely positive map and bi-optimal by optimal and define maximal completely map majorized by $\phi$ in analogy with the definition for decomposable maps, we obtain the following result.

\begin{thm}\label{thm7}
Let A be a finite dimensional $C^*$-algebra and $H$ a finite dimensional Hilbert space.  Let $\phi\colon A\to B(H)$ be a positive map.  Then there  are  a maximal completely positive map $\alpha\colon A\to B(H)$ majorized by $\phi$ and an optimal map $\beta\colon A\to B(H)$ such that $\phi = \alpha + \beta$.
\end{thm}

\section{Spin factors}

In the present section we illustrate the decomposition theorems, Theorem 5 and Theorem 7, by the projection of $B(H)$ onto a spin factor.  Following \cite{HS} we recall that a \textit spin system in $B(H) $ is a set of symmetries, i.e. self-adjoint unitaries $s_1,...,s_m$  satisfying the anti-commutation relations $s_is_j + s_js_i = 0$ for $iÊ\neq j$. Let
$$
\sigma_1 = \left(\begin{array}{cc}
1&0\\
0&-1
\end{array}\right), \quad
\sigma_2 = \left(\begin{array}{cc}
0&1\\
1&0
\end{array}\right),\quad
\sigma_3 = \left(\begin{array}{cc}
0&-i\\
i&0
\end{array}\right)
$$
denote the Pauli matrices in $M_2$.  Then we can construct a spin system $\{s_1,...,s_{2n}\}$ in $ M_{2^n} = \bigotimes_1{^n} M_2$ as follows, where $1 \leq k <  n-1$.
\begin{eqnarray*}
&&s_1 = \sigma_1 \otimes 1^{\otimes n-1}\\
&&s_2 = \sigma_2 \otimes 1^{\otimes n-1}\\
&&.\\ 
&&.\\
&&s_{2k+1}= \sigma_3^{\otimes k} \otimes \sigma_1 \otimes1^{\otimes n-k-1}\\
&&s_{2k+2}= \sigma_3^{\otimes k} \otimes \sigma_2 \otimes1^{\otimes n-k-1}\\
&&.\\
&&.\\
&&s_{2n-1} = \sigma_3^{\otimes n-1}Ê\otimes \sigma_1\\
&&s_{2n} = \sigma_3^{\otimes n-1}Ê\otimes \sigma_2
\end{eqnarray*}
where for $a \in M_2$, $a^{\otimes k} $ denotes the k-fold tensor product of $a$ with itself.

Let $V_m$ denote the linear span of $s_0= 1, s_1,...,s_m$.  Then $V_m$ is a spin factor of dimension m+1 in $M_{2^n}$.  For m=2n the $C^*$-algebra $C^*(V_m)$ generated by $V_{2n}$ equals $M_{2^n}$, so in that case $V_m$ is irreducible, see \cite{HS}, Theorem 6.2.2.  If m=2n-1 then $C^*(V_m) = M_{2^{n-1}}\bigoplus M_{2^{n-1}} \subset M_{2^n}$. 

By \cite{ES} or \cite{S} , Proposition 2.2.10, if $Tr$ denotes the usual trace on $M_{2^n}$ then there exists a positive idempotent map $P \colon M_{2^n} \to V_m + iV_m$ given by $Tr(P(a)b) = Tr(ab)$ for all $aÊ\in M_{2^n}, b \in V_m + iV_m,\ \  m\leq 2n$.  Then $P$ restricted to the self-adjoint part of $M_{2^n}$ is a projection map onto $V_m$.  With the Hilbert-Schmidt structure the set $\{1,s_1,...,s_m\}$ is an orthonormal basis for $V_m$ with respect to the normalized trace $2^{-n} Tr$ on $M_{2^n}$.  Thus $P$ has the form 
$$
P(a)  = 2^{-n} \sum_0 ^m Tr(s_i a)s_i.
$$
By \cite{R2} or \cite{S},Theorem 2.3.4 $P$ is atomic if $n \neq 2,3,5$.  By \cite{HS},Theorem 6.2.3 $V_m$ is a JW-factor of type $I_2$, i.e. for each minimal projection $e\in V_m, 1-e$  
is also a minimal projection.  Thus $Tr(e)=2^{n-1}$. Note that for each $e_i, i\geq 1, e_{+}= 1/2 (1+s_i)$ and $e_{-} = 1/2 (1 - s_i )$ are such projections.

Let $t$ denote the transpose on $M_2$.  Then 
$$
\sigma_1^t = \sigma_1,\ \  \sigma_2^t =\sigma_2, \ \ \sigma_3^t = -\sigma_3.
$$
Since the transpose on $M_{2^n}$ is the tensor  product $t^{\otimes n}$,  it follows from the defining equations for $s_k$ that
$$
s_{2k+1}^t = (-1)^k s_{2k+1},\ \   s_{2k+2}^t = (-1)^k s_{2k+2}.
$$
It follows in particular that $P \circ t = t \circ P$.

\begin{lem}\label{lem6}
Define a symmetry $W  \in M_{2^n} $ as follows:
\enu{i} If n is odd, n = 2m+1,  $W = (1 \otimes \sigma_3)^{\otimes m} \otimes 1$.
\enu{ii} If n is even, n = 2m,  $W = (1 \otimes \sigma_3)^{\otimes m}.$

Then $AdW(s_k) = s_{k}^t $ for all $1 \leq k \leq 2n$.  Hence $AdW(a) = a^t$ for all $aÊ\in V_n$. Furthermore, if n is of the form n=4m+i, i= 0,1, then $W Ê\in C^*(V_n)$.
\end{lem}
\bp
If $k=1,2$, then $AdW(s_k) = s_k = s_{k}{^t} $, so we may assume $k\geq 3$. We first consider the case when $ k=2j+1$ with j odd. Then
$$
s_k = \sigma_{3}^{\otimes j} \otimes \sigma_1 \otimes 1^{\otimes {n-j-1}}.
$$
Thus by definition of $W$, since $ Ad\sigma_3(\sigma_1) = - \sigma_1$, we have
$$
AdW(s_k) = \sigma_3^{\otimes j} \otimes (-\sigma_1) \otimes 1^{\otimes{n-j-1}} = -s_k = (-1)^j s_k = s_{k}^t.
$$
Similarly if $k=2j+2$ with j odd, then $AdW(s_k) = s_{k}{^t} $.  Now let $k=2j+1$ with j even. Then
$$
AdW(s_k) = \sigma_3^{\otimes j} \otimes \sigma_1 \otimes 1^{\otimes {n-j-1}} = s_k =  (-1)^j s_k = s_{k}^t.
$$
Similarly for $k= 2j+2$ with j even.  Thus in every case $AdW(s_k) = s_{k}^t.$ Since $V_n$ is the real linear span of $s_k$, k= 0,1,...,n, $AdW(a) = a^t$ for all $a\in V_n$

If $n=4m+i, i  =0,1$, then, since $\sigma_3^{t} = -\sigma_3$, and there are $2m$ factors of $\sigma_3$ in $W$,  we have $W^t = W$.  If$ i=0$ then by \cite{HS}, Theorem  6.2.2, $C^*(V_n) = M_{2^n}$, so clearly $W  \in C^*(V_n)$.  If $n=4m+1$ then again by \cite{HS}, Theorem 6 2.2
$$
C^*(V_{4m+1}) = M_{2^{4m}} \bigoplus M_{2^{4m}} \subset M_{2^{4m+1}}.
$$
Since in this case
$ W = (1 \otimes  { \sigma_3 })^{\otimes 2m } \otimes 1$, 
it follows that
$W  \in M_{4^m} \otimes {\0C}  \subset C^*(V_{4m+1}) =C^*(V_n),$
completing the proof of the lemma.

\begin{lem}\label{lem7}
Let $m\leq 2n$ and $P \colon M_{2^n} \to V_m$ be the trace invariant projection.  Let $W$ be as in Lemma 7. Then
$$
P = P \circ t \circ AdW.
$$
\end {lem}
\bp
By Lemma 8 if $aÊ\in V_m$ then $tÊ\circ AdW(a) = a.$  Thus if $x \in M_{2^n}$ then
$$
(P \circ t \circ AdW)  \circ (P \circ t \circ AdW)(x)  = P \circ (P \circ t \circ AdW)(x) = P \circ t \circ AdW(x).
$$
Thus $P \circ t \circ AdW$ is idempotent with range $V_m$ and being the identity on $V_m$.  Since $P$ is trace invariant, if $x \in M_{2^n}, y \in V_m$ we have
\begin{eqnarray*}
Tr(P \circ t \circ AdW(x) y) &=& Tr( t \circ AdW(x) y)  = Tr(AdW(x) y^t)\\
& &= Tr(x AdW \circ t(y)) = Tr(xy)  = Tr(P(x)y),
\end{eqnarray*}
using that $AdW \circ t = t \circ AdW = \iota$ on $V_m$, where $\iota$ is the identity map on $V_m$.  The lemma follows.

\medskip  

The following lemma is probably well known, but is included for completeness.

\begin{lem}\label{lem8}
Let $a \in B(H)$ be positive and $e,f$ projections in $B(H)$ with sum 1.  Then
$$
2(eae + faf ) \geq a.
$$
\end{lem}
\bp
We have
$$
a = (e + f)a(e + f) = eae + eaf + fae + faf.
$$
Let
$$
b = (e - f)a(e - f)  = eae - eaf - fae + faf \geq 0.
$$
Thus 
$$
a \leq a + b = 2(eae + faf),
$$
as asserted.

We shall need the following slight extension of a result of Robertson \cite{R}.  For simplicity we show it in the finite dimensional case.  Recall that $M^\prime$ denotes the commutant  for a set $M \subset B(H)$ and that $B_{sa}$ denotes the set of self-adjoint operators in $M$.

\begin{lem}\label{lem9}
Let $H$ be a finite dimensional Hilbert space, and let $B Ê\subset  B(H)$ be a $C^*$-algebra and $ A \subset B_{sa}$ a Jordan algebra with $1 \in A$.  Suppose $P  \colon B_{sa} \to  A$ is a positive projection  map.  Suppose $\phi  \leq P$ is a completely positive map, $\phi \colon B  \to B$. Then $\phi(1) \in C^*(A)^\prime$, and $\phi(x) = \phi(1)x$ for  $x \in C^*(A)$.
\end{lem}
\bp
By \cite{S}, Lemma 2.3.5, since $P(x) = x$ for $x \in A, \phi(1) \in A$, and $\phi(x) = \phi(1)x = x\phi(1),$ for $x\in A$.  Since $C^*(A)$ is the $C^*$-algebra generated by $A, \phi(1) \in C^*(A)^\prime$.  Since $H$ is finite dimensional, if $e$ is the range projection of $\phi(1), \phi(1)$ has a bounded inverse $\phi(1)^{-1}$ on $eH$.  Thus 
$$
\psi = \phi(1)^{-1} e \phi 
$$
is a unital map of $B$ into $eBe$ such that for $x \in A$,
$$
\psi(x) = \phi(1)^{-1} e \phi(x) = \phi(1)^{-1} \phi(1) x = ex.
$$
Thus  $\psi_{ \vert A}$ is a Jordan homomorphism, so $A \subset D = \{x \in B_{sa}: \psi(x^2) = \psi(x)^2\} $,  the definite set for $\psi$.  Since $\psi$ is completely positive, by \cite {R} or  \cite{S}, Proposition 2.1.8, $D$ is the self-adjoint part of a $C^*$ -algebra, hence $\psi$ is a homomorphism on $C^*(A)$.  Since by the above $\psi(x) = ex$ for $x \in A,  \psi(x) = ex$ for $x \in C^*(A)$.  If $x \in C^*(A), 0 \leq x \leq 1$ then $\phi(x) \leq \phi(1) =e\phi(1)$.  Thus $\phi(x) = e\phi(x)$, so that for all $x \in C^*(A)$, we have
$$
\phi(x) = e\phi(x) = \phi(1) \psi(x) = \phi(1)ex = \phi(1) x,
$$
proving the lemma.

\begin{lem}\label{lem10}
Let $P \colon M_{2^n}  \to V_m,  m \leq 2n$ be the trace invariant projection  Then $P \geq  2^{-n} \iota$, and $P \geq 2^{-n} t \circ  AdW,$ with $W$ as in Lemma 8.  Furthermore there exists a 1-dimensional projection $q \in M_{2^n} $ such that $P(q) = 2^{-n}1$, hence
$$
2^{-n} = max \{ \lambda \geq 0: P \geq \lambda \iota\}.
$$
\end{lem}
\bp
Let $p$ be a 1-dimensional projection in $M_{2^n}$. Since $V_m$ is a JW-factor of type $I_2$, \cite{HS}, Theorem 6.1.8, there are two minimal projections $e$ and $f$ in $V_m$ with sum 1 and $a,b \geq 0$ such that 
$$
P(p) = ae + bf.
$$
By \cite{S}, Proposition 2.1.7, $ P(epe) = eP(p)e = ae$, so that 
$$
a 2^{n-1} = Tr(ae) = Tr(P(epe)) = Tr(epe).
$$
Hence
$$
a = 2^{-n+1} Tr(epe),    b = 2^{-n+1} Tr(fpf).
$$
Since $epe$ is positive of rank 1, $Tr(epe) \geq epe$.  Thus, using Lemma 10 we get
\begin{eqnarray*}
P(p) &=& 2^{-n+1} ( Tr(epe)e + Tr(fpf) f)\\
&\geq& 2^{-n+1}( epe + fpf)\\
&\geq& 2^{-n+1} \2 (epe +  epf + fpe + fpf)\\
&=& 2^{-n} p.
\end{eqnarray*}
Since this holds for all 1-dimensional projections $p, PÊ\geq 2^{-n} \iota$.  By Lemma 9 it thus follows that 
$$
P = P \circ t \circ AdW \geq 2^{-n} t \circ AdW,
$$
proving the first part of the lemma.

To show the second part we exhibit a 1-dimensional projection $q$ such that $ P(q) = 2^{-n} 1$.  The Pauli matrix $\sigma_3$ is of the form  $\sigma_3 = e_0 - f_0 \in M_2$ with $e_0, f_0$ 1-dimensional projections in $M_2$.  Let $Tr_2$ denote the usual trace on $M_2$.  Then for j=1,2, we have
\begin{eqnarray*}
0&=& Tr_2(\sigma_3 \sigma_j) = Tr_2(e_0 \sigma_j) - Tr_2(f_0 \sigma_j)\\
&=& Tr_2(e_0 \sigma_j  - (1- e_0)\sigma_j)\\
&=& 2 Tr_2(e_0 \sigma_j) - Tr_2(\sigma_j)\\
&=& 2 Tr_2(e_0 \sigma_j).
\end{eqnarray*}
Furthermore, $Tr_2(e_0 \sigma_3) = Tr_2(e_0(e_0 - f_0)) = Tr_2(e_0) = 1.$   Let $q = e_0 ^{\otimes n}  \in M_{2^n}.$  If $ j=2k-i, i=0,1$, then $s_j = \sigma_3^{\otimes k-1} \otimes \sigma_j  \otimes 1^{\otimes{n-k} }$.  From the above we thus have 
$$
Tr(q s_j) = Tr_2(e_0 \sigma_j) = 0.
$$
Thus, since $s_0 = 1$, we have
$$
P(q) = 2^{-n} (\sum_0^m Tr(q s_j)s_j) = 2^{-n} Tr(q s_0)s_0 = 2^{-n} 1,
$$
completing the proof.

\medskip

The projection $q$ above is not symmetric because $\sigma_3^{t} = -\sigma_3 = f_0 - e_0$, so that $e_0^{t} = f_0.$   Furthermore $AdW(q) = AdW( e_0^{\otimes n}) = q$, hence $t \circ AdW(q) = q^t \perp q$.   These properties of q will limit our choice of $V_m$ in our study of $P$.

In the case $m=2^n$ there are four classes of non-isomorphic irreducible Jordan subalgebras  of $(M_m)_{sa}$, namely $(M_m)_{sa} $ itself, $V_{2n}$, $S_m$, the real symmetric matrices in $M_m$, and $M_{2^{n-1}}(\0H)_{sa}$,  the self-adjoint $2^{n-1} \times 2^{n-1}$ matrices over the quaternions $\0H$ represented as $2 \times 2$ matrices, see \cite {HS}, Ch. 6.  Presently we shall specialize to the case when $V_{2n} \subset (M_{2^{n-1}})_{sa}$. We refer the reader to \cite{LM} for further information on this case. 

With our previous notation with $W$ defined as in Lemma 8 let
$$
Q(X) = \2 (x + t\circ AdW(x)).
$$
Then $Q$ is the projection of $M_{2^n}$ onto the fixed point set of the anti-automorphism $t\circ AdW$, hence by Lemma 8 is the projection onto the reversible Jordan algebra $A_{2n}$ containing $V_{2n}$.  Thus, if $V_{2n} \subset M_{2^{n-1}}(\0H)_{sa}$ then $Q\colon M_{2^n} \to M_{2^{n-1}}(\0H)_{sa}$.

\begin{lem}\label {lem 11}
With the above notation, if $V_{2n} \subset A_{2n} = M_{2^{n-1}}(\0H)_{sa}$ and $P$ the projection $P \colon M_{2^n} \to V_{2n}$, then
$$
P = P\vert_{A_{2n}}  \circ  Q  \geq 2^{-n+1} Q.
$$
\end{lem}

\bp
It suffices to show $P(p) \geq 2^{-n+1} p$ for all minimal projections $p$ in $A_{2n}$.  For such a projection $Tr(p) = 2$. We have $P(p) = ae + bf$, $ a, b \geq 0$, as in the proof of Lemma 12. Then $a= 2^{-n+1} Tr(epe), b= 2^{-n+1} Tr(fpf)$.  Since $p$ is a minimal projection in $A_{2n}$, $pep = \lambda p, pfp = \mu p$ with $\lambda, \mu \geq 0.$ Then 
$$
(epe)^2 = epepe = \lambda epe.
$$
Since rank $epe$ = rank $  pep = 2$, $epe = \lambda_0   e_0$ with $e_0$ a projection in $A_{2n}$ of dimension 2.  Thus
$$
\lambda_0 2 = Tr(\lambda_0 e_0) = Tr(epe) = Tr(pep) = Tr(\lambda p) = \lambda 2.
$$
Therefore $\lambda_0 = \lambda$.  Thus $epe = \lambda e_0$, and similarly $fpf = \mu f_0$.  We thus have, since $e\geq e_0$ and $f\geq f_0$,
\begin{eqnarray*}
P(p) &=& 2^{-n+1} (Tr(epe)e  + Tr(fpf)f) \\
&=& 2^{-n+1}(Tr(\lambda e_0)e + Tr(\mu f_0)f)\\
&\geq& (2\lambda e_0 + 2\mu f_0)\\
&=& 2^{-n+1}(2 epe +2 fpf)\\
&\geq& 2^{-n+1}(epe + epf + fpe + fpf)\\
&=& 2^{-n+1} p,
\end {eqnarray*}
where we used Lemma 10. The proof is complete.

\begin{lem}\label{lem13}
Given $V_{2n}$ and $A_{2n}$ as above, and assume $A_{2n} \cong M_{2^{n-1}}(\0H)_{sa}$.  Then there exists a 1-dimensional projection $q$ in $M_{2^n}$ such that $Q(q) = \2(q + q^t)$ with $q \bot q^t$, $P(q) = 2^{-n}1$, and $\beta = P - 2^{-n+1} Q $ is bi-optimal.
\end{lem}
\bp
By Lemma 13 $P\vert_{A_{2n}} \geq 2^{-n+1}\iota$. Since $P = P \circ Q$ we therefore have $\beta = P\circ Q - 2^{-n+1} Q \geq 0$.  $V_{2n}$ is irreducible by \cite{HS}, Theorem 6.2.2, so $C^*(V_{2n}) = M_{2^n}$, so by Lemma 12 there is a 1-dimensional $q\in C^*(V_{2n})$ such that $2^{-n} 1 = P(q) = P(Q(q))$.  By the comments after Lemma 12, $q^t = t\circ AdW(q) \bot q$, so in particular 
$$
Q(q) = \2(q + \circ AdW(q) ) = \2(q + q^t).
$$
Furthermore
$$
\beta(Q(q)) = P(Q(q)) - 2^{-n+1} Q(q) = 2^{-n} (1 - (q + q^t)).
$$
To show $\beta$ is bi-optimal, let $\phi \leq \beta$ be completely positive.  Then by Lemma 11, $\phi(x) = \phi(1)x = \lambda x, \lambda \geq 0$, since $\phi(1) \in C^*(V_{2n})^\prime = \0C$. Thus
$$
\lambda(q+q^t) = \phi(q+q^t) = 2\phi(Q(q)) \leq 2\beta(Q(q)) = 2^{-n} (1 - (q+q^t)).
$$
Since $q+q^t \bot 1-(q+q^t), \lambda = 0$, so $\phi = 0$. Thus $\beta$ is optimal.

Next, if $\phi\leq \beta$ is copositive, then $t\circ \phi$ is completely positive,and 
$$
t\circ \phi \leq  t\circ P = P \circ t = P \circ AdW,
$$
since $P = P \circ t \circ AdW$ by Lemma 9. Thus by Lemma 11 $t\circ \phi = \lambda \iota$ with $\lambda \geq 0$. Hence
$$
\lambda(q+q^t) = t\circ \phi(q+q^t) = 2 t\circ \phi(Q(q)) \leq2 t\circ \beta(Q(q)) = 2^{-n}(1-(q+q^t))^t =2^{-n} (1 - (q+q^t)),
$$
so again $\lambda =0$, and $\phi=0$.  Thus $\beta$ is bi-optimal, completing the proof to the lemma.

\medskip
From the above we see that if $\phi \leq P$ is completely positive or copositive, then $\phi \leq \lambda Q$ for some $\lambda \geq 0$.  Since $P \geq \alpha = 2^{-n+1} Q$, and $P(q) = 2^{-n} 1$, it follows that $\alpha$ is a maximal decomposable map majorized by $P$.

Summarizing Lemma 14 and the above comments we obtain the following result.

\begin{thm}\label{thm 14}
Assume the reversible Jordan algebra $A_{2n}$ containing $V_{2n}$ is isomorphic to $M_{2^{n-1}}(\0H)_{sa}$, and let $Q\colon M_{2^n} \to A_{2n}$ be the trace invariant projection. Let $\alpha = 2^{-n+1} Q$ and $\beta = P - \alpha$.  Then $P = \alpha + \beta$ is a decomposition as in Theorem 5.
\end{thm}

The following result describes Theorem 7 in detail for $P$.

\begin{thm}\label{17}
Let $P\colon M_{2^n} \to V_{2n}$ be the trace invariant projection.  Let $\alpha = 2^{-n} \iota$, and $\beta = P - 2^{-n} \iota$, where $\iota$ is the identity map.  Then $\alpha$ is a maximal  completely positive map majorized by $P$,  $\beta$ is optimal, and $P = \alpha + \beta$.
\end{thm}

\bp
By Lemma 12 $P \geq \alpha$, so $\beta \geq 0$, and there exists a 1-dimensional projection $q\in M_{2^n}$ such that $P(q)= 2^{-n}1$.  Since $V_{2n}$ is irreducible the argument in the proof of Lemma 14 shows that if $\phi \leq \beta$ is completely positive, then $\phi = \lambda \iota$ with $\lambda \geq 0$. Thus
$$
\lambda q = \phi(q) \leq \beta(q) = 2^{-n} 1 - 2^{-n}q = 2^{-n}(1 - q),
$$
which implies $\lambda =0$.  Thus $\beta$ is optimal.  As remarked before the statement of Theorem 15 $\alpha$ is a maximal completely positive map majorized by $P$.  The proof is complete.

It was crucial in the proof of Theorem 15 that $A_{2n} = M_{2^{n-1}}(\0H)_{sa}$, so $dim q =2$ for a minimal projection $q$ in $A_{2n}$.  In the case when $A_{2n} = S_{2^n}$, the real $2^n \times 2^n$ matrices, we have been unable to find a 1-dimensional projection $p\in A_{2n}$ such that $P(p)=2^{-n} 1$, so that for each minimal projection $e\in V_{2n}$ we have
$$
Tr(pe) = Tr(epe) = Tr(P(epe)) = Tr(eP(p)e) = Tr(e2^{-n}1) = \2,
$$
so $Tr(p . )$ is the trace on $V_{2n}$.

If $n=1,  V_2 = S_2 =A_1$, so $Tr(p . )$ is never a trace on $A_1$.  We next show this for $V_4$ too, showing in particular the well known result that $A_2 = M_2(\0H)_{sa}$.  We thus leave it as an open question whether there is n  such that $Tr(p . )$ can be a trace on $V_{2n}$ for a 1-dimensional projection $p\in A_{2n}$, or even for $p\in M_{2^n}$.

\begin {ex}.
If n=2 then there is no positive rank 1 operator  $x \in M_4$ such that $t \circ AdW(x) = x$.
\end{ex}

\bp
  Let $\bar{\phi}  \colon M_2 \to M_2$ be defined by
$$
\bar{\phi} \left(\begin{array}{cc}
a & b\\
c & d
\end{array}\right) 
= \left(\begin{array}{cc}
d & -c\\
-b & a
\end{array}\right)
$$
Then $\bar{\phi} = Ad\sigma_3$ as is easily seen. Let $\phi = t  \circ  {\bar{\phi}}. $ Then $\phi$ is an anti-automorphism of order 2, and 
$$
\phi \left(\begin{array}{cc}
a & b\\
c & d
\end{array}\right) 
= \left(\begin{array}{cc}
d & -b\\
-c & a
\end{array}\right)
$$
is such that  $\1R = \{A \in M_2: \phi(A^*) = A\} $ is the quaternions.  Also $\phi = t \circ \sigma_3$.  
For simplicity of notation let $\rho = Ad\sigma_3$.  Let $T$ denote the $4 \times 4 $ matrix
$$
\left(\begin{array}{cc}
A & B\\
C & D
\end{array}\right)
$$
with $A,B,C,D \in M_2$. Then 
$$
\iota \otimes \rho(T^*)  = 
\left(\begin{array}{cc}
\rho(A)^*  & \rho(C)^*\\
\rho(B)^*  & \rho(D)^*
\end{array}\right)
$$
Therefore 
$$
t \circ (\iota \otimes \rho) (T^*) = 
\left(\begin{array}{cc}
t \circ \rho(A)^* & t \circ \rho(B)^*\\
t \circ \rho(C)^* & t \circ \rho(D)^*
\end{array}\right)
$$
Thus $t \circ (\iota \otimes \rho) (T^*) = T$ if and only if
$$
A = \phi(A^*), B = \phi(B^*), C= \rho(C^*), D = \phi(D^*)
$$
if and only if $A,B,C,D \in \0H$,  and so $T \in M_2{(\0H)}.$ But $M_2{(\0H)}$ contains no positive rank 1 operators, so there is no positive rank 1 $x \in M_4$ such that $t \ \circ AdW(x) =x$, completing the proof of the example.

 \medskip
 
 If $\1P =\{ s_i : i  \in \0N \}$ is an infinite spin system then the norm closed linear span $V_{\infty}$ of 1 and $\1P$ is the infinite spin factor.  The $C^*$-algebra $C^*(V_{\infty})$ generated by $V_{\infty}$ is the CAR-algebra $A$ which is isomorphic to the infinite tensor product of $M_2$ with itself, see e.g. \cite {HS}, Theorem 6.2.2.  By \cite {ES}, Lemma 2.3, there exist a unique trace invariant positive projection $P$ of $C^*(V_{\infty})_{sa}$ onto $V_{\infty}$.  If $M_{2^n} = \otimes_1^n M_2$ is imbedded in $C^*(V_{\infty})$ by $x \rightarrow x \otimes 1 \in M_{2^n} \otimes \otimes_{n+1}^{\infty} M_2$ , it is clear that $P\vert _{M_{2^n}} = P_n,$ the trace invariant projection onto $V_{2n}$.  Thus if $\phi \leq P$ is decomposable then $\phi\vert_{M_{2^n}} \leq P\vert_{M_{2^n}} = P_n$ for n even. Thus by Lemmas11 and 12, $\phi\vert_{M_{2^n}} \leq 2^{-n} \iota\vert_{M_{2^n}}$.  But if $m\geq n$ is even then 
 $$
 \phi\vert_{M_{2^n}} = (\phi\vert_{M_{2^m}})\vert_{M_{2^n}} \leq 2^{-m} (\iota\vert_{M_{2^m}})\vert_{M_{2^n}}.
 $$
Thus 
$$ 
 \phi\vert_{M_{2^n}}  \leq 2^{-m} \iota\vert_{M_{2^n}}  
 $$ for all even $m \geq n$.  Thus $\phi= 0$.   Similarly if $\phi \leq t \leq P$.  We have thus shown
 
 \begin{cor}\label{cor14}
 Let $P$ be the projection of the self-adjoint part of the CAR-algebra onto the spin factor $V_\infty$.  Then $P$ is bi-optimal.
 \end{cor}

Department of Mathematics,  University of Oslo, 0316 Oslo, Norway.

e-mail  erlings@math.uio.no

\end{document}